\documentclass[11pt]{amsart}

\usepackage{amsmath,amssymb,amsthm,amsrefs,a4wide}
\usepackage{latexsym}
\usepackage{indentfirst}
\usepackage{graphicx}
\usepackage{placeins}
\usepackage{booktabs}
\usepackage{algorithm}
\usepackage{algorithmic}
\usepackage{multirow}

\theoremstyle{plain}

\theoremstyle{definition}

\theoremstyle{remark}

\DeclareMathOperator{\argmin}{arg min}

\renewcommand{\d}{\mathrm{d}}
\newcommand{\KL}{\mathrm{KL}}
\newcommand{\rKL}{\mathrm{rKL}}
\newcommand{\Hlg}{\mathrm{H}}

\newcommand{\diag}{\mathrm{diag}}

\newcommand{\bbm}{\begin{bmatrix}}
\newcommand{\ebm}{\end{bmatrix}}

\newcommand{\R}{\mathrm{R}}

\newcommand{\tp}{\tilde{p}}
\newcommand{\tg}{\tilde{g}}

\begin{document}

\title[Mirror Descent for Interacting Free Energy]{Mirror Descent Algorithms for Minimizing \\Interacting Free Energy}

\author[]{Lexing Ying}
\address[Lexing Ying]{Department of Mathematics and ICME, Stanford University, Stanford, CA 94305}
\email{lexing@stanford.edu}

\thanks{The work of L.Y. is partially supported by the U.S. Department of Energy, Office of Science,
  Office of Advanced Scientific Computing Research, Scientific Discovery through Advanced Computing
  (SciDAC) program and also by the National Science Foundation under award DMS-1818449. The author
  thanks Wuchen Li and Wotao Yin for comments and suggestions.}


%\subjclass[2010]{65Z05, 82B28, 82B80.}

\begin{abstract}
  This note considers the problem of minimizing interacting free energy. Motivated by the mirror
  descent algorithm, for a given interacting free energy, we propose a descent dynamics with a novel
  metric that takes into consideration the reference measure and the interacting term. This metric
  naturally suggests a monotone reparameterization of the probability measure. By discretizing the
  reparameterized descent dynamics with the explicit Euler method, we arrive at a new
  mirror-descent-type algorithm for minimizing interacting free energy. Numerical results are
  included to demonstrate the efficiency of the proposed algorithms.
\end{abstract}

\maketitle

%-----------------------------------
\section{Introduction}\label{sec:intro}

This paper considers the problem of minimizing free energies of the following form
\begin{equation} \label{eqn:F}
  F(p) = D(p||\mu) + \int_\Omega p(x) V(x) \d x + \frac{1}{2}\iint p(x) W(x,y) p(y) \d x \d y
\end{equation}
for a probability density $p$ over domain $\Omega$. $D(p||\mu)$ is a divergence function between $p$
and a reference density $\mu$ and typically examples are Kullback-Leibler divergence, reverse
Kullback-Leibler divergence, and Hellinger divergence. In the interacting term $\iint p W p \d x \d
y$, $W$ is symmetric and can either be positive-definite or not. Non-positive-definite interacting
terms appear in Keller-Segel models in mathematical biology and granular flows in kinetic theory.
Recently, positive-definite interacting terms appear in the mean field modeling of neural network
training \cite{chizat2018global,mei2018mean,rotskoff2018neural,sirignano2018mean}.

%goal
The goal of this paper is to develop fast first-order algorithms for identifying minimums of
\eqref{eqn:F}. When $F$ is convex (for example, when $W$ is positive-definite), there exists a
unique global minimizer and the goal is to compute this global minimizer efficiently. When $F$ is
non-convex, there are typically many local minimums and the more moderate goal is to find one such
local minimum.

%difficulty
There are several difficulties for computing local minima for \eqref{eqn:F}. First, this is an
optimization problem over probability simplex, hence one needs to deal with the constraints $p(x)\ge
0$ and $\int p(x) \d x=1$. Second, when the reference measure $\mu(x)$ varies drastically for
different $x\in\Omega$, the optimization problem can be quite ill-conditioned. Third, we aim to
avoid costly second-order Newton or quasi-Newton methods that involve matrix inversions or solves.

%motivation
\subsection{Motivations and approach}
Our approach is motivated by the mirror descent algorithm \cite{nemirovsky1983problem} popularized
recently in the machine learning community. Because of several nice computational and analytical
features, the mirror descent algorithm has played a significant role in online learning and
optimization. For an objective function $E(p)$ over the space of probability densities, it finds a
minimizer of $E(p)$ as follows. Given a current density $p^k$, each step solves for
\begin{equation}\label{eqn:mdprox}
\tp = \argmin_p E(p^k) + \frac{\delta E}{\delta p}(p^k) \cdot(p-p^k) + \frac{1}{\eta} D_{\KL}(p||p^k)
\end{equation}
and then projects $\tp$ back to the space of probability densities. Taking derivative of
\eqref{eqn:mdprox} in $p$ results in
\[
\eta \frac{\delta E}{\delta p}(p^k) + \ln(\tp/p^k) + 1 = 0,
\]
with $\tp$ proportional to $p^k \exp\left(-\eta \frac{\delta E}{\delta p}(p^k)\right)$.  Projecting
it back to the probability simplex via rescaling gives
\begin{equation}\label{eqn:md}
  p^{k+1} =\frac{1}{Z} p^k \exp\left(-\eta \frac{\delta E}{\delta p}(p^k)\right),\quad
  Z = \int p^k \exp\left(-\eta \frac{\delta E}{\delta p}(p^k)\right) \d x.
\end{equation}

Let us now give a different derivation of the mirror descent algorithm from a more numerical
analysis perspective. The starting point is the natural gradient flow of $E(p)$ with the Fisher-Rao
metric $\diag(1/p)$:
\[
\dot{p} = - \frac{1}{1/p} \left(\frac{\delta E}{\delta p} +c\right) = -p\left(\frac{\delta E}{\delta p} +c\right),
\]
where $\frac{\delta E}{\delta p}$ is Frechet derivative and $c$ is the Lagrange multiplier
associated with $\int_\Omega p(x) \d x=1$. Moving $p$ to the left hand side gives rise to
an equation of $\ln p$.
\[
\dot{(\ln p)} = -\left(\frac{\delta E}{\delta p} +c\right).
\]
Using the explicit Euler method in the new variable $\ln p$ with step size $\eta$ results in
\[
\ln p^{k+1} = \ln p^k -\eta \left(\frac{\delta E}{\delta p}(p^k) +c\right),
\]
where $c$ is determined from the condition $\int p^{k+1} \d x = 1$ and this is equivalent to
\eqref{eqn:md}. This derivation shows that mirror descent can be viewed as the explicit Euler
discretization of the natural gradient flow in the reparameterization $\phi(p) \equiv \ln p$.

The mirror descent is effective when the Hessian of the energy function $E(p)$ is close to the
Fisher-Rao metric $1/p$, up to a constant scaling. This is the case for
\[
E(p) = \int p(x)\ln p(x) \d x + \int V(x) p(x) \d x,
\]
where the Hessian is exactly the Fisher-Rao metric. In this case, the natural gradient is
\[
\dot{(\ln p)} = - (\ln p + V +c). %\quad \dot g = - (g+V+c).
\]
This is a linear system of ordinary differential equations with coefficient $1$ in the new variable
$\ln p$. The stiffness is gone and one can take large steps.

%In fact, in this case, taking step size equal to one is equivalent to the Newton method in the
%original $p$ variable.

%current approach
Coming back to the free energy \eqref{eqn:F}, the mirror descent algorithm described above is not
particularly effective, due to the existence of the reference measure $\mu$ (in the reverse KL and
Hellinger cases) as well as the extra interacting term $W$. In fact, for general $\mu$ and $W$, the
Fisher-Rao metric $1/p$ in the natural gradient algorithm is quite far away from the Hessian matrix
of the Newton method. Therefore, there is no reason to expect the standard mirror descent algorithm
to be efficient. Our approach consists of the following steps:
\begin{itemize}
\item Choose an appropriate diagonal metric based on $\mu$ and $W$;
\item Design a reparameterization function $\phi$ based on the chosen metric;
\item Derive the algorithm by performing the explicit Euler discretization;
\item Work out the renormalization step.
\end{itemize}

\subsection{Related work.}
The mirror descent algorithm \cite{nemirovsky1983problem,beck2003mirror} was proposed as an
effective first-order method for convex optimization by taking into consideration the geometry of
the problem. For certain types of constraint sets, the mirror descent algorithm is nearly optimal
among first order methods \cite{bubeck2015convex}, offering an almost dimensional independent
convergence rate. In the setting of online optimization, mirror descent also allows one to obtain a
bound for the cumulative regret \cite{arora2012multiplicative,bubeck2011introduction}. There is a
vast literature on mirror descent and related algorithms and we refer to
\cite{bubeck2015convex,shalev2012online} for further discussions.

The interacting free energy of form \eqref{eqn:F} appear in several applications, such as
Keller-Segel models \cite{perthame2006transport} in mathematical biology, as well as the granular
flow in kinetic theory \cite{carrillo2003kinetic,villani2006mathematics}. In these applications, the
evolution of the probability density is governed by the Wasserstein gradient flow
\cite{jordan1998variational,otto2001geometry} of the free energy, i.e., the gradient flow with
respect to the Wasserstein metric $-\nabla\cdot(p\nabla (\cdot))$. The main computational task in
these applications is to compute the evolution of the Wasserstein gradient flow and several
numerical methods based on finite element, finite volume, and particle methods
\cite{bessemoulin2012finite,carrillo2019blob,li2018natural,li2019fisher,liu2018positivity} have been
proposed for this.  Compared with these algorithms, the goal of this paper is different as we only
care about the minimizers. Therefore, we have the freedom to pick any descent dynamics that leads to
the minimizer. As we have seen, our flow is closer to the natural gradient rather than the
Wasserstein gradient.

\subsection{Contents.}
The paper considers three common cases of the divergence term $D(p||\mu)$ and is organized as
follows. Section \ref{sec:kl} addresses the Kullback-Leibler divergence, Section \ref{sec:rkl} is
about the reverse Kullback-Leibler case, and finally Section \ref{sec:hlg} discusses the Hellinger
distance case. In each case, we address both the case of positive-definite $W$ term as well as the
general situation of non-positive-definite $W$.

As the metric adopted here is of the Fisher-Rao type as opposed to the Wasserstein type, there is no
derivative involved in the computation. To simplify the presentation and also to make connection
with the numerical implementation, we work with a probability density $\{p_1,\ldots,p_n\}$ over a
discrete set of $n$ points $\{x_1,\ldots, x_n\}$ rather than over a continuous space. The
interacting free energy can be written as
\[
F(p) = D(p||\mu) + \sum_i p_i V_i + \frac{1}{2}\sum_{ij} p_i W_{ij} p_j.
\]
This is indeed the setup when \eqref{eqn:F} is discretized with a numerical treatment.

%-----------------------------------
\section{Kullback-Leibler divergence} \label{sec:kl}

For the KL divergence case,
\[
D_{\KL}(p||\mu) = \sum_{i=1}^n p_i \ln p_i/\mu_i = \sum_{i=1}^n p_i \ln p_i -\sum_{i=1}^n p_i \ln \mu_i.
\]
The second term can be absorbed into the potential $V$ and hence it is equivalent to consider
\[
F_{\KL}(p) = \sum_i p_i \ln p_i + \sum_i V_i p_i + \frac{1}{2}\sum_{i,j} p_i W_{ij} p_j.
\]
The Hessian is given by
\[
\frac{\delta^2 F_{\KL}}{\delta p^2} = \diag\left(\frac{1}{p}\right) + W.
\]
When $W$ is non-positive-definite, the safe way is to just use $\diag(1/p)$ as the gradient
metric. When $W$ is positive-definite, we extract the diagonal $\alpha = \diag(W)\in \R^n$ of $W$
and use $\diag(1/p+\alpha)$ as the gradient metric.

%-------
\subsection{Non-positive-definite case}\label{sec:kl_non}

Using $\diag(1/p)$ as the metric, the gradient flow is
\[
\dot{p} = -p (\ln p + V + Wp +c).
\]
Moving the metric to the left hand side gives
\[
\dot{(\ln p)} = -(\ln p + V + Wp +c).
\]
If we introduce a reparameterization from $p\in\R^n$ to $g\in\R^n$ with $g_i=\phi_i(p_i)\equiv\ln
p_i$ and $p_i=\phi_i^{-1}(g_i) = \exp(g_i)$
\begin{align*}
  & \phi_i: p_i\rightarrow g_i,\quad (0,1)\rightarrow (-\infty, 0), \\
  & \phi_i^{-1}: g_i\rightarrow p_i,\quad (-\infty,0)\rightarrow (0,1),
\end{align*}
the gradient flow becomes
\[
\dot g = -(g + V + W p +c).
\]
The explicit Euler discretization gives
\begin{align*}
  & \tg = g^k -\Delta t (g^k + V + Wp^k),\\
  & g^{k+1} = \tg+c.
\end{align*}
The constant $c$ is determined by the normalization condition
\[
\sum_i \phi_i^{-1}(\tg_i+c) = 1,
\]
which leads to $c = -\ln \left(\sum_i \exp(\tg_i)\right).$

\begin{figure}[h!]
  \centering
  \includegraphics[width=0.32\textwidth]{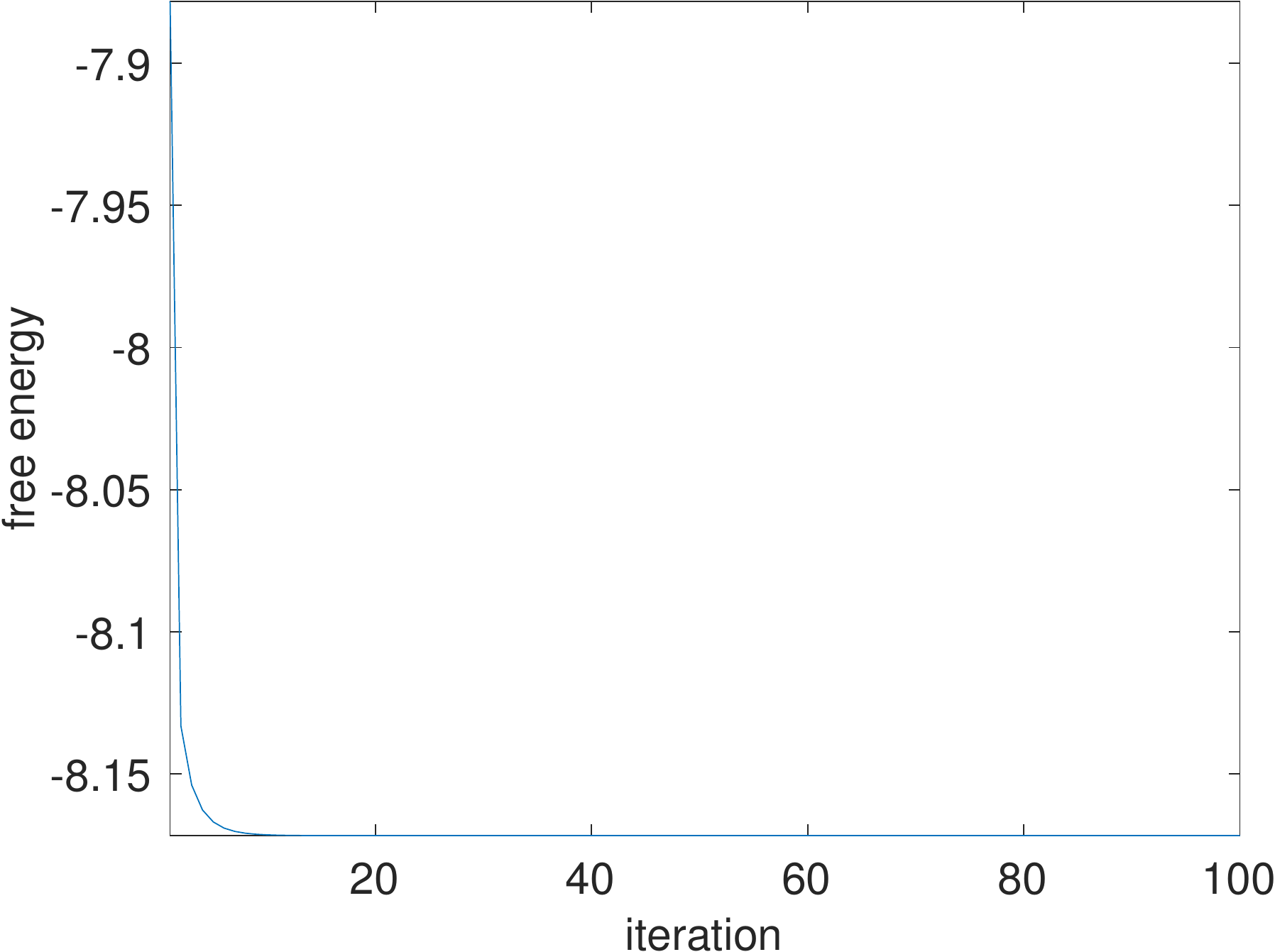}
  \includegraphics[width=0.32\textwidth]{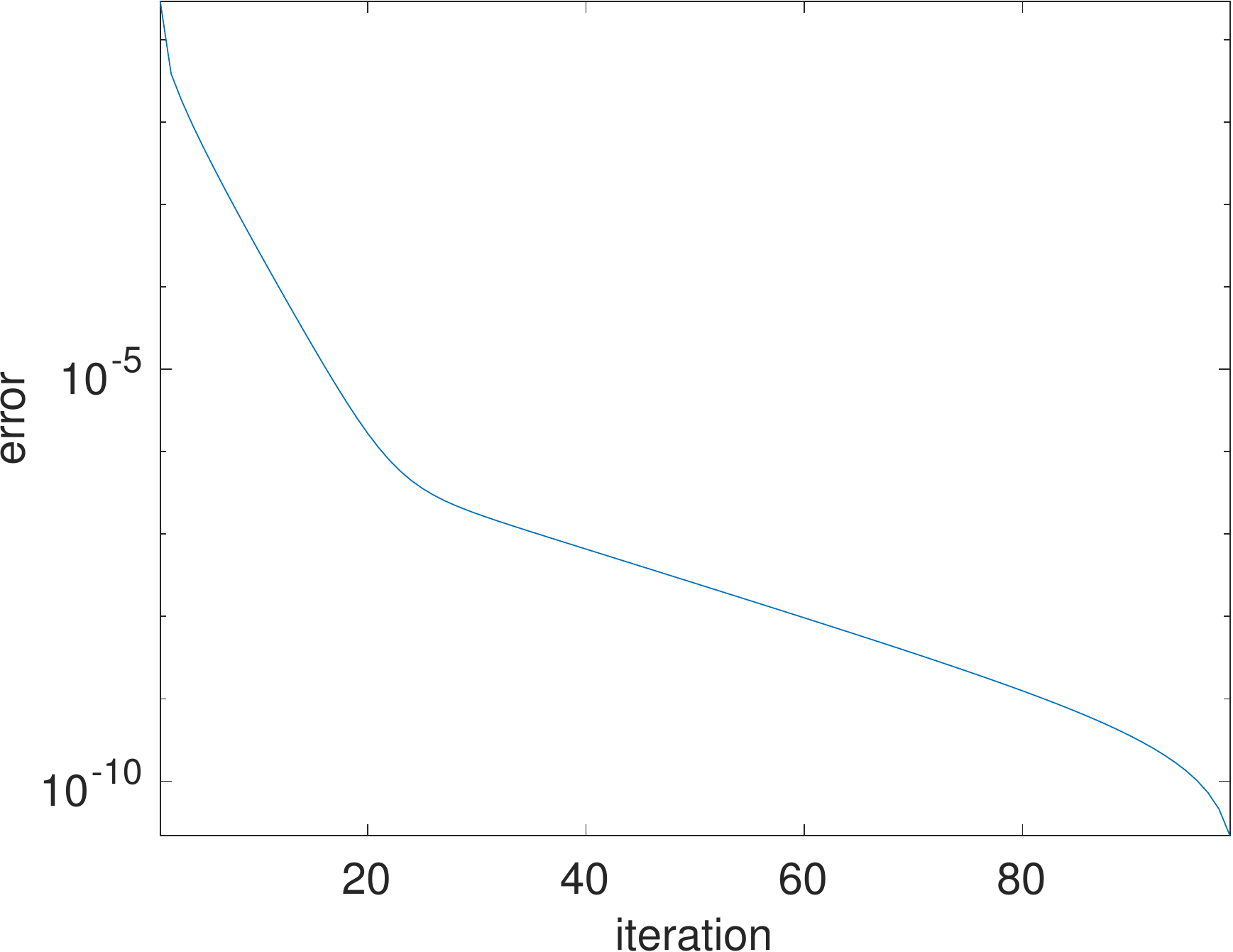}
  \includegraphics[width=0.32\textwidth]{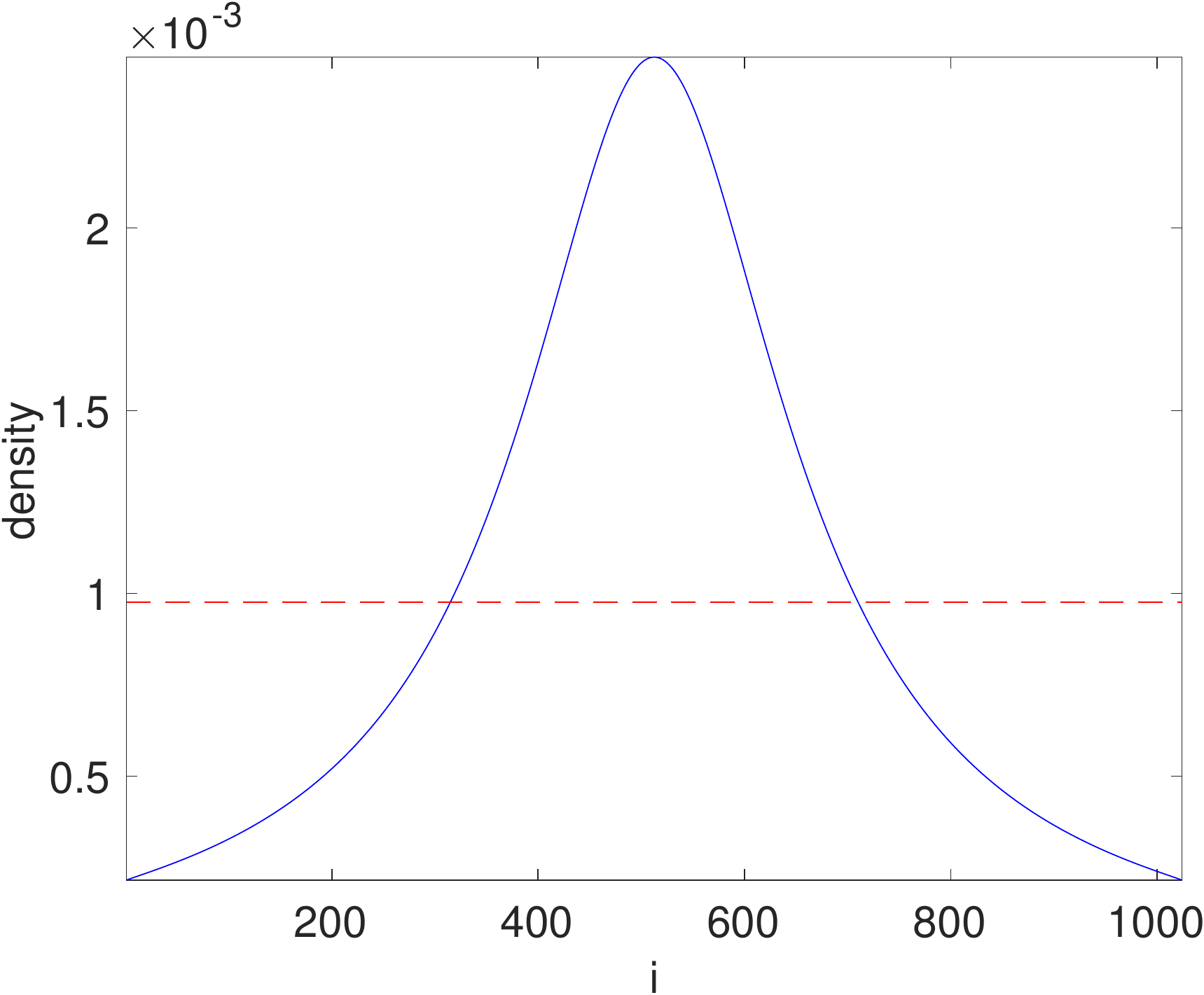}
  \caption{KL divergence, non-positive-definite case with a Keller-Segel free energy.  Left: free
    energy vs. iteration. Middle: free energy error vs. iteration. Right: density $p$ at the final
    iteration (solid line) and the uniform density (dashed line). The uniform density is the
    minimizer when $W$ term is absent.}
  \label{fig:kl_non}
\end{figure}

We illustrate the efficiency of the algorithm with a Keller-Segel model.  Consider the domain
$[0,1]$ discretized with $n=1024$ points $\{x_i=\frac{i}{n}\}$. The potential $V$ is zero and the
interacting term is
\[
W_{ij} = \frac{3}{2} \ln(|x_i-x_j|+\varepsilon)
\]
with $\varepsilon=10^{-6}$. The step size $\Delta t$ is taken to be $1$. Starting from a random
initial condition, we run the descent algorithm for $100$ steps. The results are summarized in
Figure \ref{fig:kl_non}. At the end of the $100$ iterations, the error is of order $10^{-10}$. The
final density shows the concentration property of the Keller-Segel free energy.

%-------
\subsection{Positive-definite case}\label{sec:kl_spd}

Using $\diag(1/p)+\alpha$ as the metric, the gradient flow is
\[
\dot{p} = -\frac{1}{1/p + \alpha} (\ln p + V + Wp +c).
\]
Moving the metric to the left hand side gives
\[
\dot{(\ln p +\alpha p)} = -(\ln p + \alpha p + V + (W-\alpha)p +c).
\]
If we introduce a reparameterization from $p\in\R^n$ to $g\in\R^n$ with $g_i = \phi_i(p_i)
\equiv\ln(p_i) + \alpha_i p_i$
\begin{align*}
  & \phi_i: p_i\rightarrow g_i,\quad (0,1)\rightarrow (-\infty, \alpha_i),\\
  & \phi_i^{-1}: g_i\rightarrow p_i,\quad (-\infty,\alpha_i)\rightarrow (0,1),
\end{align*}
the gradient flow becomes
\[
\dot g = -(g + V + (W-\alpha) p +c).
\]
The explicit Euler discretization gives
\begin{align*}
  &\tg = g^k -\Delta t (g^k + V + (W-\alpha)p^k),\\
  & g^{k+1} = \tg+c.
\end{align*}
The constant $c$ is determined by the normalization condition
\[
\sum_i \phi_i^{-1}(\tg_i+c) = 1.
\]
Let us observe that $\sum_i \phi_i^{-1}(\tg_i+c)$ is an increasing function in $c$ as each
$\phi_i^{-1}$ is increasing. The correct value $c$ can be shown to be in
\[
\left( \min\left(\ln \frac{1}{n} + \frac{\alpha_i}{n} - \tg_i\right), \min(\alpha_i - \tg_i) \right).
\]
Plugging the two endpoints of the interval shows that at the left endpoint
$\sum_i\phi_i^{-1}(\tg_i+c)<1$ and at the right endpoint $\sum_i \phi_i^{-1}(\tg_i+c)>1$. Therefore,
there is a unique $c$ value satisfies $\sum_i \phi_i^{-1}(\tg_i+c) = 1$ within this interval. This
can be easily found using Newton, bisection, or interpolation methods \cite{forsythe1977computer}.

\begin{figure}[h!]
  \centering
  \includegraphics[width=0.32\textwidth]{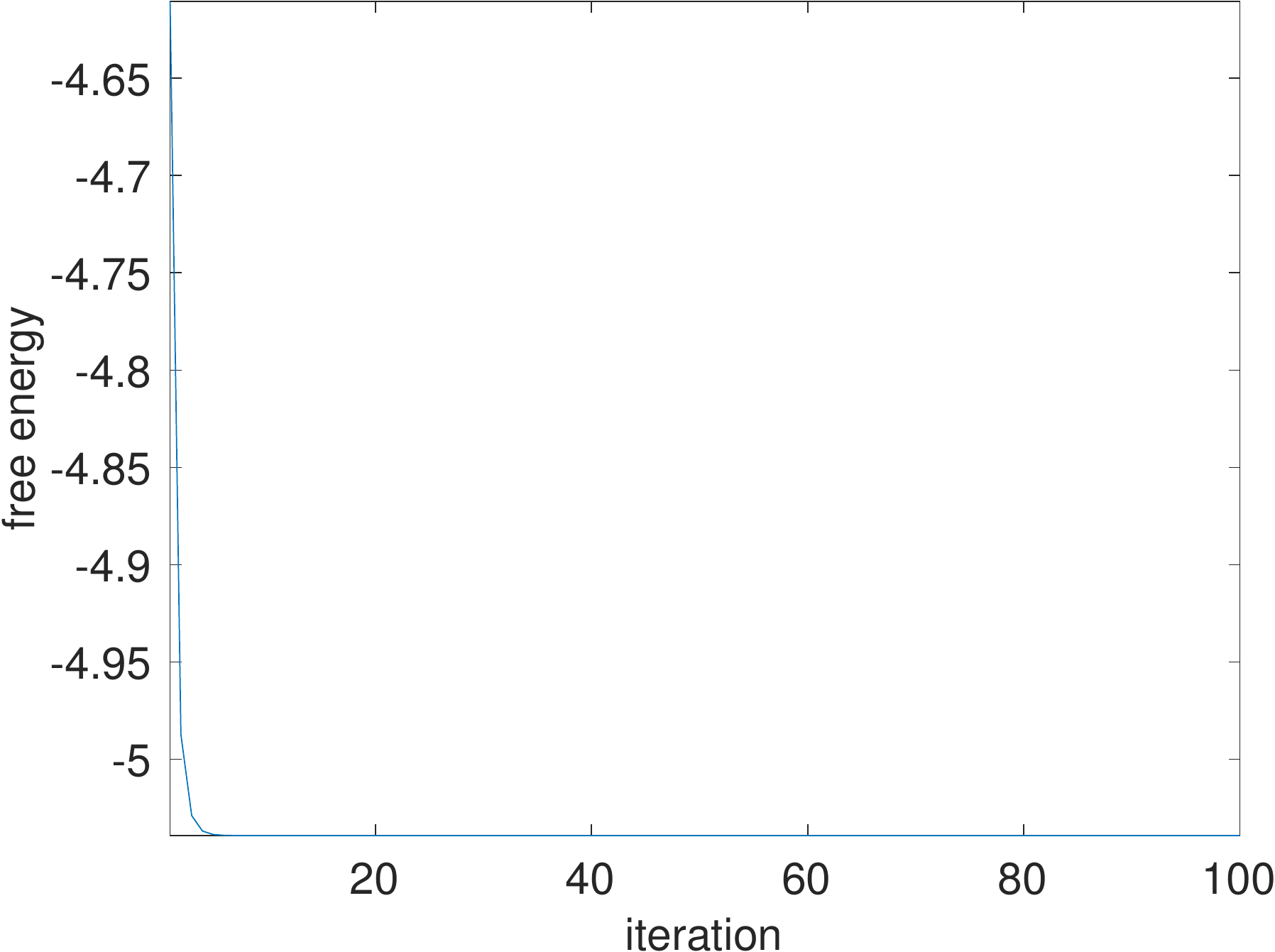}
  \includegraphics[width=0.32\textwidth]{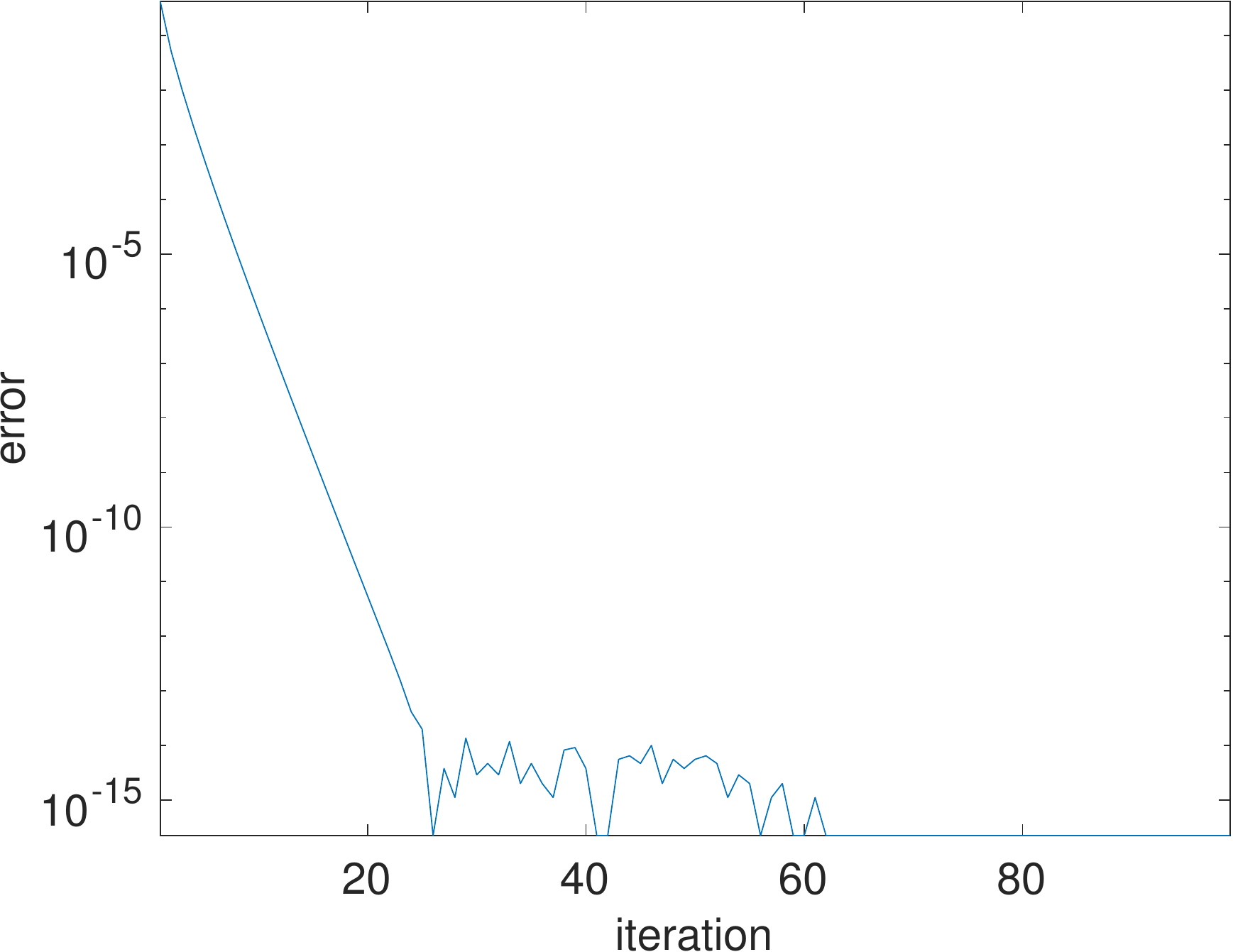}
  \includegraphics[width=0.32\textwidth]{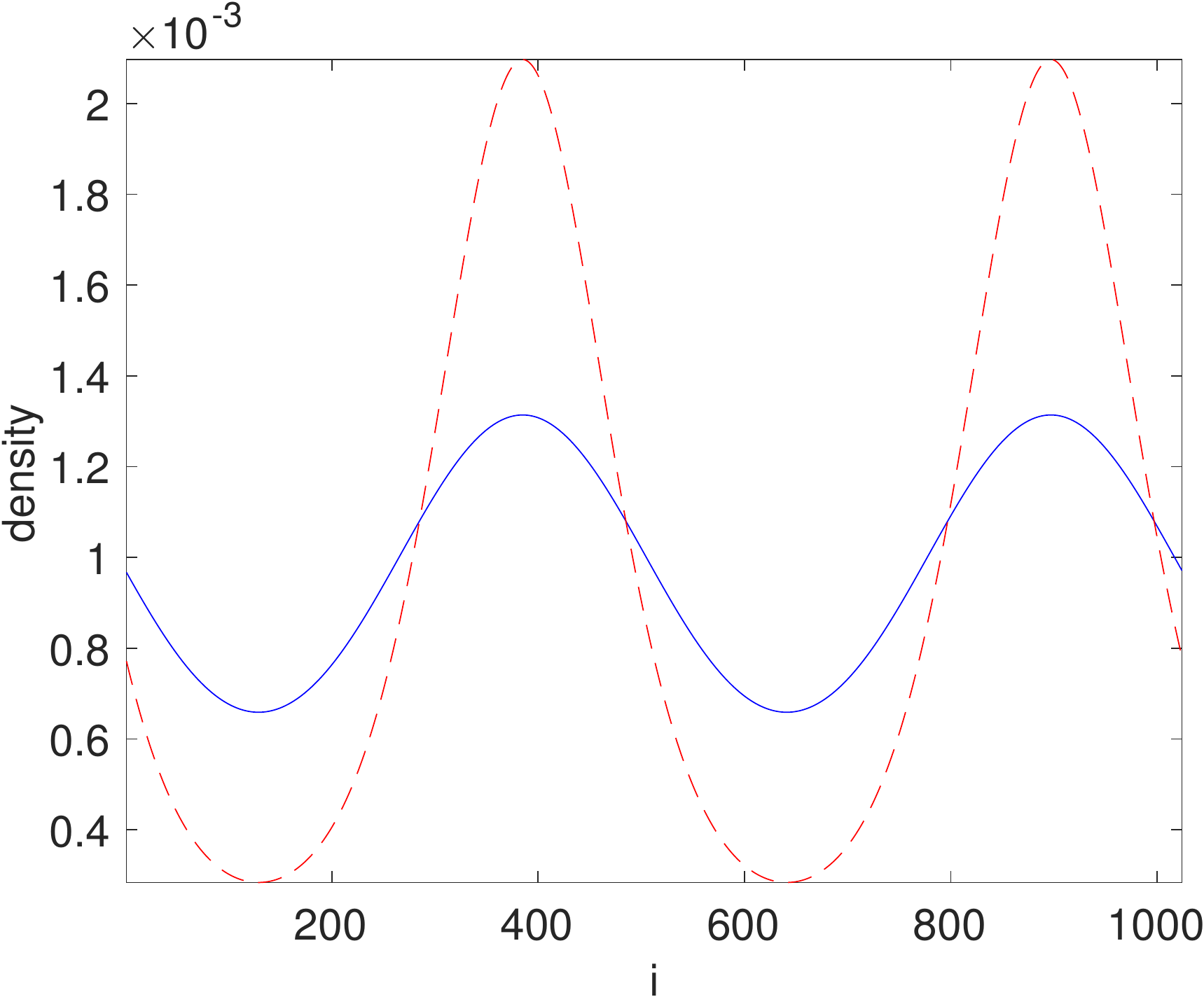}
  \caption{KL divergence, positive-definite case.  Left: free energy vs. iteration. Middle: free
    energy error vs. iteration. Right: density $p$ at the final iteration (solid line) and the
    minimizer density with $W=0$.}
  \label{fig:kl_spd}
\end{figure}

To illustrate the efficiency of this algorithm, we consider the periodic domain $[0,1]$ discretized
with $n=1024$ points. The potential $V$ is chosen to be $V_i=\sin(4\pi x_i)$ and the interacting
term is
\[
  W_{ij} =
  \begin{cases}
    \alpha, & i=j,\\
    \alpha/2, & i=j\pm 1,\\
    0, & \text{otherwise},
  \end{cases}
\]
with $\alpha=10^3$. Hence, $\alpha_i=10^3$ for each $i$. The step size $\Delta t$ is taken to be
$1$. Starting from a random initial condition, we run the algorithm for $100$ steps with the results
summarized in Figure \ref{fig:kl_spd}. Within $20$ iterations, it reaches within $10^{-15}$
accuracy. The final probability density shows that the interacting term in the free energy further
suppresses oscillations in the minimizing density.

%-----------------------------------
\section{Reverse Kullback-Leibler divergence} \label{sec:rkl}

For the reverse KL divergence
\[
D_{\rKL}(p||\mu) = \sum_i \mu_i \ln \mu_i/p_i = \sum_i \mu_i \ln\mu_i - \sum_i \mu_i \ln p_i.
\]
The free energy is now
\[
F_\rKL(p) = -\sum_i \mu_i \ln p_i + \sum_i V_i p_i + \frac{1}{2}\sum_{i,j} p_i W_{ij} p_j.
\]
The Hessian is given by
\[
\frac{\delta^2 F_{\rKL}}{\delta p^2} = \diag\left( \frac{\mu}{p^2}\right) + W
\]
and it can be quite far from the mirror descent choice $\diag\left(1/p^2\right)$ even when $W$ is
zero, since $\mu$ can be drastically different for different $i$. When $W$ is non-positive-definite,
it is safe to continue using $\diag\left(\mu/p^2\right)$ as the gradient metric. When $W$ is
positive-definite, we extract the diagonal $\alpha = \diag(W)$ of $W$ and use
$\diag\left(\mu/p^2 + \alpha\right)$ as the gradient metric.

%-------
\subsection{Non-positive-definite case}\label{sec:rkl_non}

Using $\diag\left(\mu/p^2\right)$ as the metric, the gradient flow is
\[
\dot{p} = -\frac{1}{\mu/p^2} \left(-\frac{\mu}{p} + V + Wp +c \right).
\]
Moving the metric to the left hand side gives
\[
\dot{(-\mu/p)} = -(-\mu/p + V + Wp +c).
\]
If we introduce a reparameterization from $p\in\R^n$ to $g\in\R^n$ with
$g_i=\phi_i(p_i)\equiv-\mu_i/p_i$ and $p_i=\phi_i^{-1}(g_i) = -\mu_i/g_i$
\begin{align*}
  & \phi_i: p_i\rightarrow g_i,\quad (0,1)\rightarrow (-\infty, -\mu_i),\\
  & \phi_i^{-1}: g_i\rightarrow p_i,\quad (-\infty,-\mu_i)\rightarrow (0,1),
\end{align*}
the gradient flow becomes
\[
\dot g = -(g + V + W p +c).
\]
The explicit Euler discretization gives
\begin{align*}
  &\tg = g^k -\Delta t (g^k + V + Wp^k),\\
  & g^{k+1} = \tg+c.
\end{align*}
The constant $c$ is determined by the normalization condition
\[
\sum_i \phi_i^{-1}(\tg_i+c) = 1,
\]
Since each $\phi_i^{-1}$ is increasing, $\sum_i \phi_i^{-1}(\tg_i+c)$ is an increasing function in
$c$. We claim that the correct value $c$ can be shown to be in
\[
\left( \min\left(-\tg_i-n\mu_i\right), \min(-\tg_i-\mu_i) \right).
\]
Plugging the two endpoints of the interval shows that at the left endpoint
$\sum_i\phi_i^{-1}(\tg_i+c)<1$ and at the right endpoint $\sum_i \phi_i^{-1}(\tg_i+c)>1$. Therefore,
there is a unique $c$ value satisfies $\sum_i \phi_i^{-1}(\tg_i+c) = 1$ within this interval.

\begin{figure}[h!]
  \centering
  \includegraphics[width=0.32\textwidth]{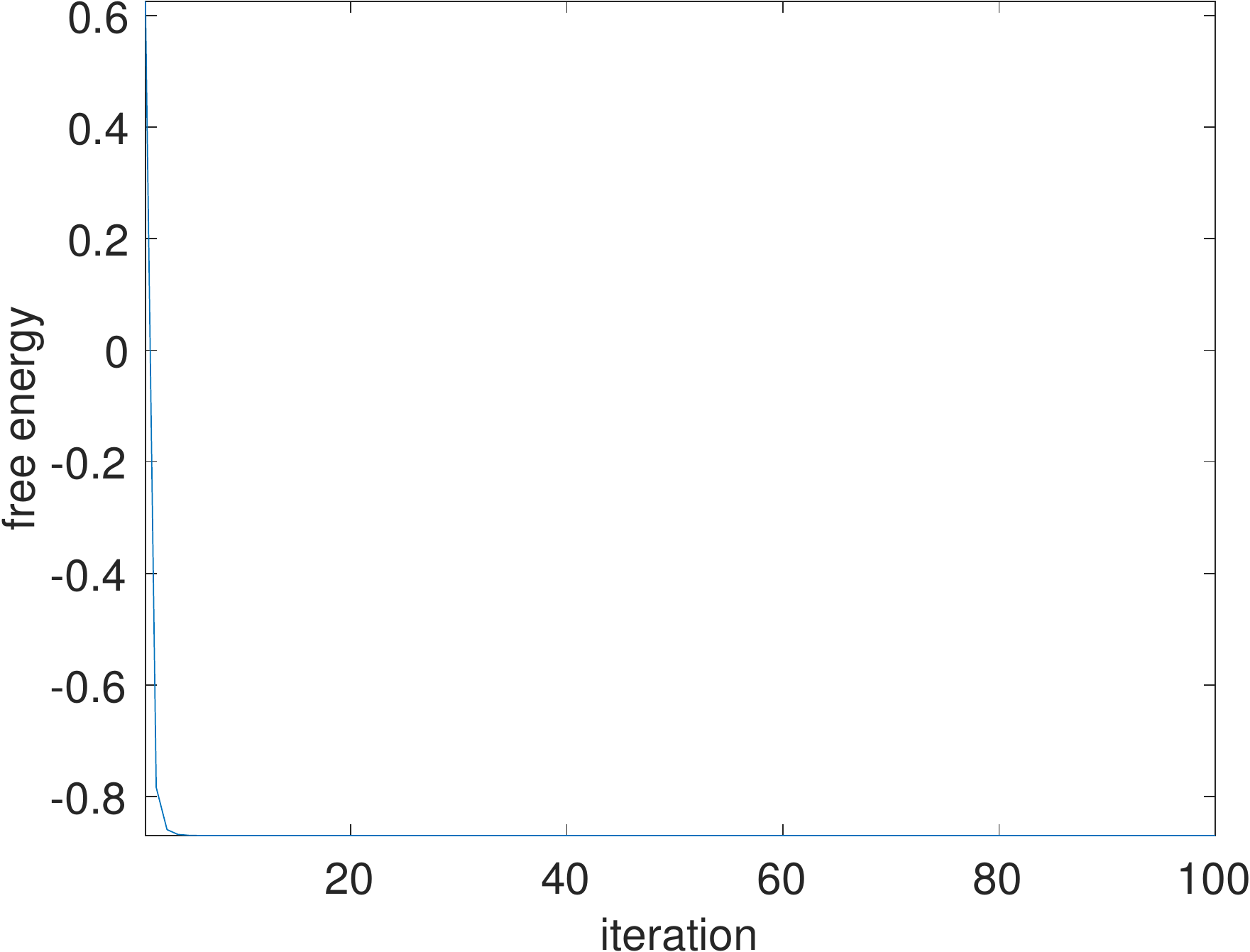}
  \includegraphics[width=0.32\textwidth]{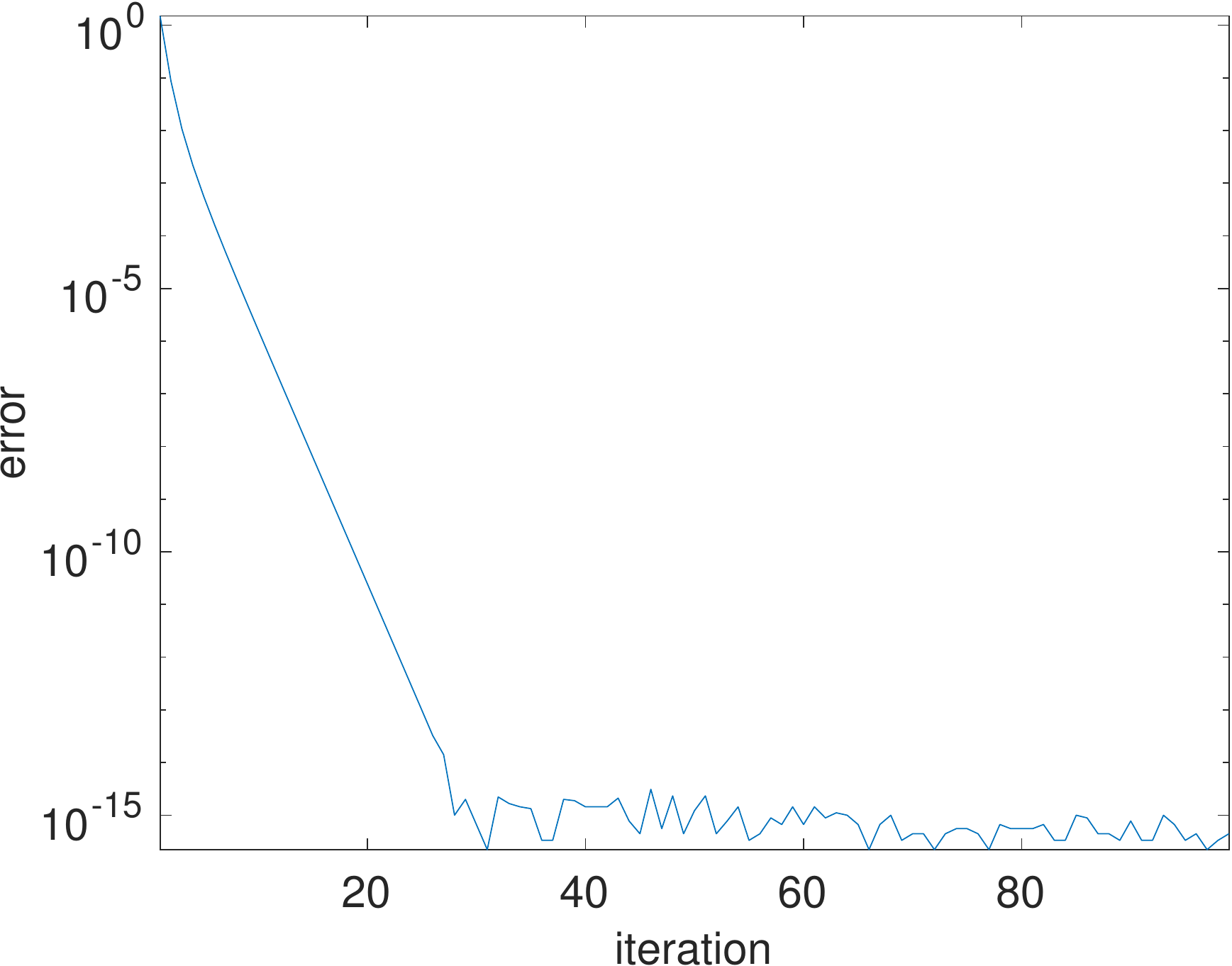}
  \includegraphics[width=0.32\textwidth]{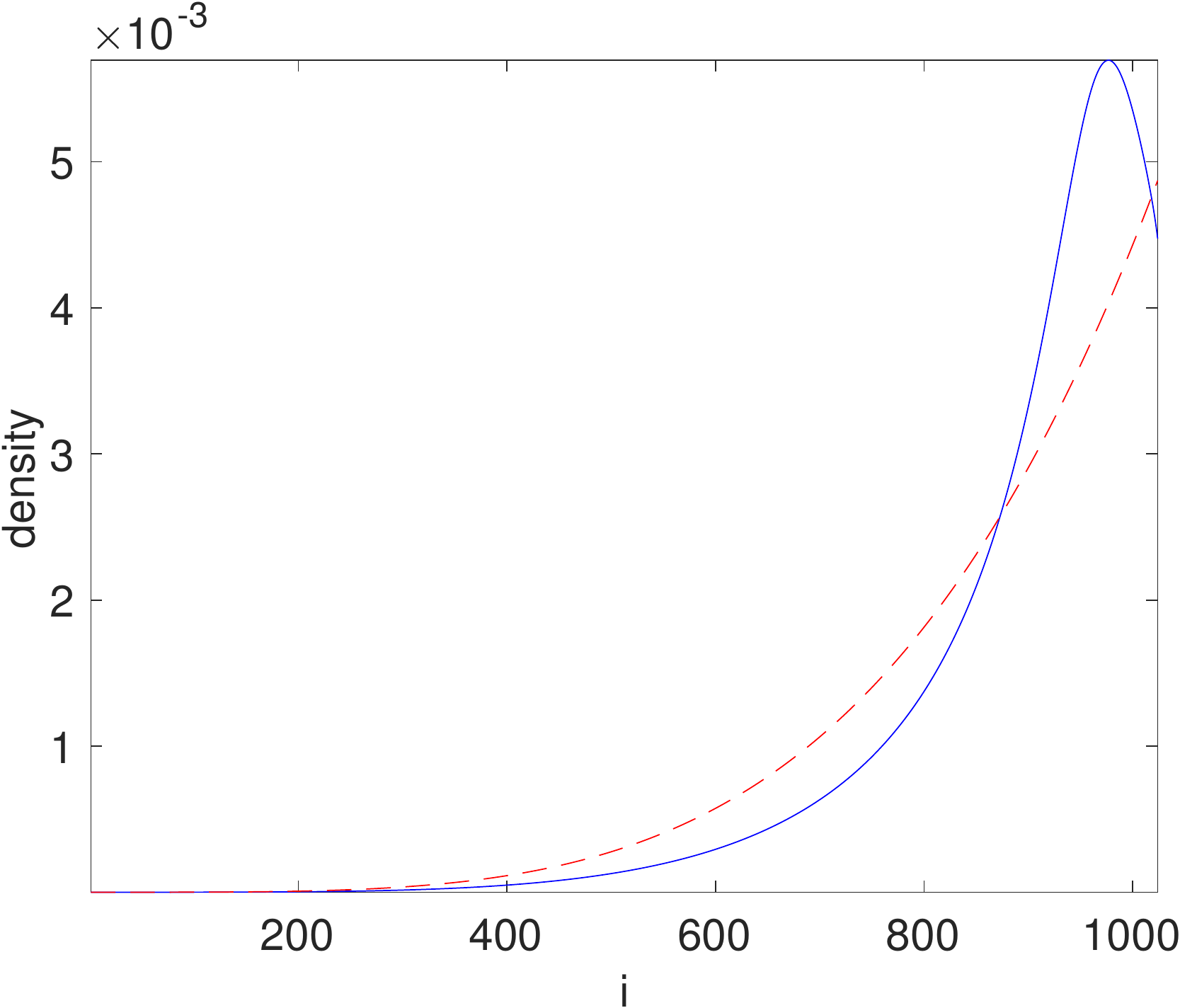}
  \caption{Reverse KL divergence, non-positive-definite case with a Keller-Segel free energy.  Left:
    free energy vs. iteration. Middle: free energy error vs. iteration. Right: density $p$ at the
    final iteration (solid line) and the reference measure $\mu$ (dashed line). The reference
    density is the minimizer when $W=0$.}
  \label{fig:rkl_non}
\end{figure}

As a numerical example, we consider a Keller-Segel model. Consider the domain $[0,1]$ discretized
with $n=1024$ points $\{x_i=\frac{i}{n}\}$. The potential $V$ is equal to zero and the interacting
term $W_{ij}$ is given by
\[
W_{ij} = \frac{2}{3} \ln(|x_i-x_j|+\varepsilon)
\]
with $\varepsilon=10^{-6}$. The reference measure $\mu_i$ is taken to be
\[
\mu_i \sim x_i^4,
\]
leading to a ratio of $10^{12}$ between the largest and the smallest $\mu_i$ values. The step size
$\Delta t$ is taken to be $1$. Starting from a random initial condition, we run the descent
algorithm for $100$ steps and the results are summarized in Figure \ref{fig:rkl_non}. Within $30$
iterations, the algorithm reaches within $10^{-15}$ accuracy.

%-------
\subsection{Positive-definite case}\label{sec:rkl_spd}

Using $\diag\left(\mu/p^2+\alpha\right)$ as the metric, the gradient flow is
\[
\dot{p} = -\frac{1}{\mu/p^2 + \alpha} (\ln p + V + Wp +c).
\]
Moving the metric to the left hand side gives
\[
\dot{(-\mu/p +\alpha p)} = -(-\mu/p + \alpha p + V + (W-\alpha)p +c).
\]
If we introduce a reparameterization from $p\in\R^n$ to $g\in\R^n$ with
$g_i=\phi_i(p_i)\equiv-\mu_i/p_i+\alpha_i p_i$ and
$p_i=\phi_i^{-1}(g_i)=\frac{g_i+\sqrt{g_i^2+4\alpha_i\mu_i}}{2\alpha_i}$
\begin{align*}
  & \phi_i: p_i\rightarrow g_i,\quad (0,1)\rightarrow (-\infty, -\mu_i+\alpha_i),\\
  & \phi_i^{-1}: g_i\rightarrow p_i,\quad (-\infty,-\mu_i+\alpha_i)\rightarrow (0,1),
\end{align*}
the gradient flow becomes
\[
\dot g = -(g + V + (W-\alpha) p +c).
\]
The Explicit Euler discretization gives
\begin{align*}
  &\tg = g^k -\Delta t (g^k + V + (W-\alpha)p^k),\\
  & g^{k+1} = \tg+c.
\end{align*}
The constant $c$ is determined by the normalization condition
\[
\sum_i \phi_i^{-1}(\tg_i+c) = 1,
\]
which can be solved since it is monotone. The correct value $c$ can be shown to be in
\[
\left( \min\left(-\tg_i-n\mu_i+\frac{\alpha_i}{n}\right), \min(-\tg_i-\mu_i+\alpha_i) \right).
\]
Plugging the two endpoints of the interval shows that the left endpoint
$\sum_i\phi_i^{-1}(\tg_i+c)<1$ and at the right endpoint $\sum_i \phi_i^{-1}(\tg_i+c)>1$. Therefore,
there is a unique $c$ value satisfies $\sum_i \phi_i^{-1}(\tg_i+c) = 1$ within this interval. 

\begin{figure}[h!]
  \centering
  \includegraphics[width=0.32\textwidth]{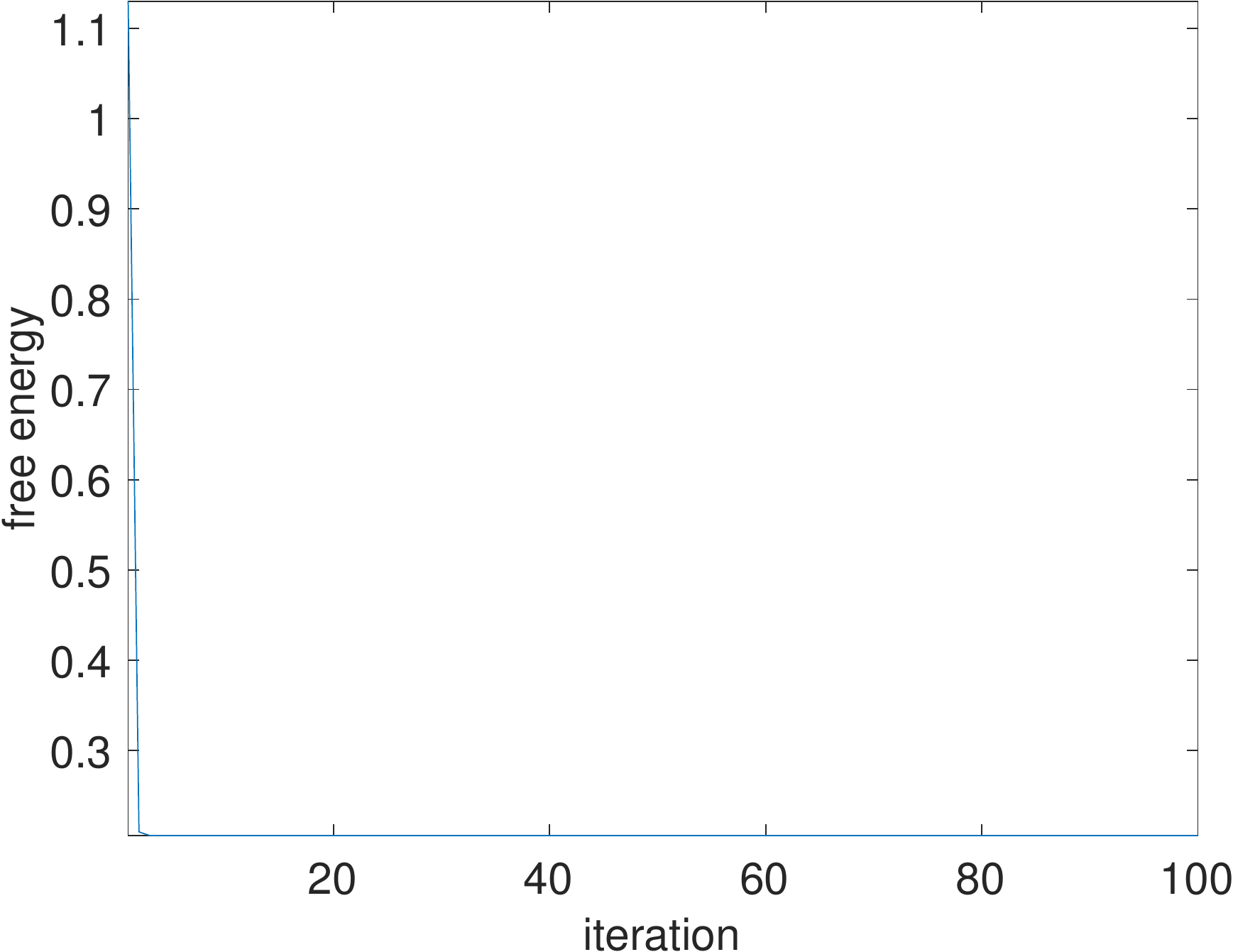}
  \includegraphics[width=0.32\textwidth]{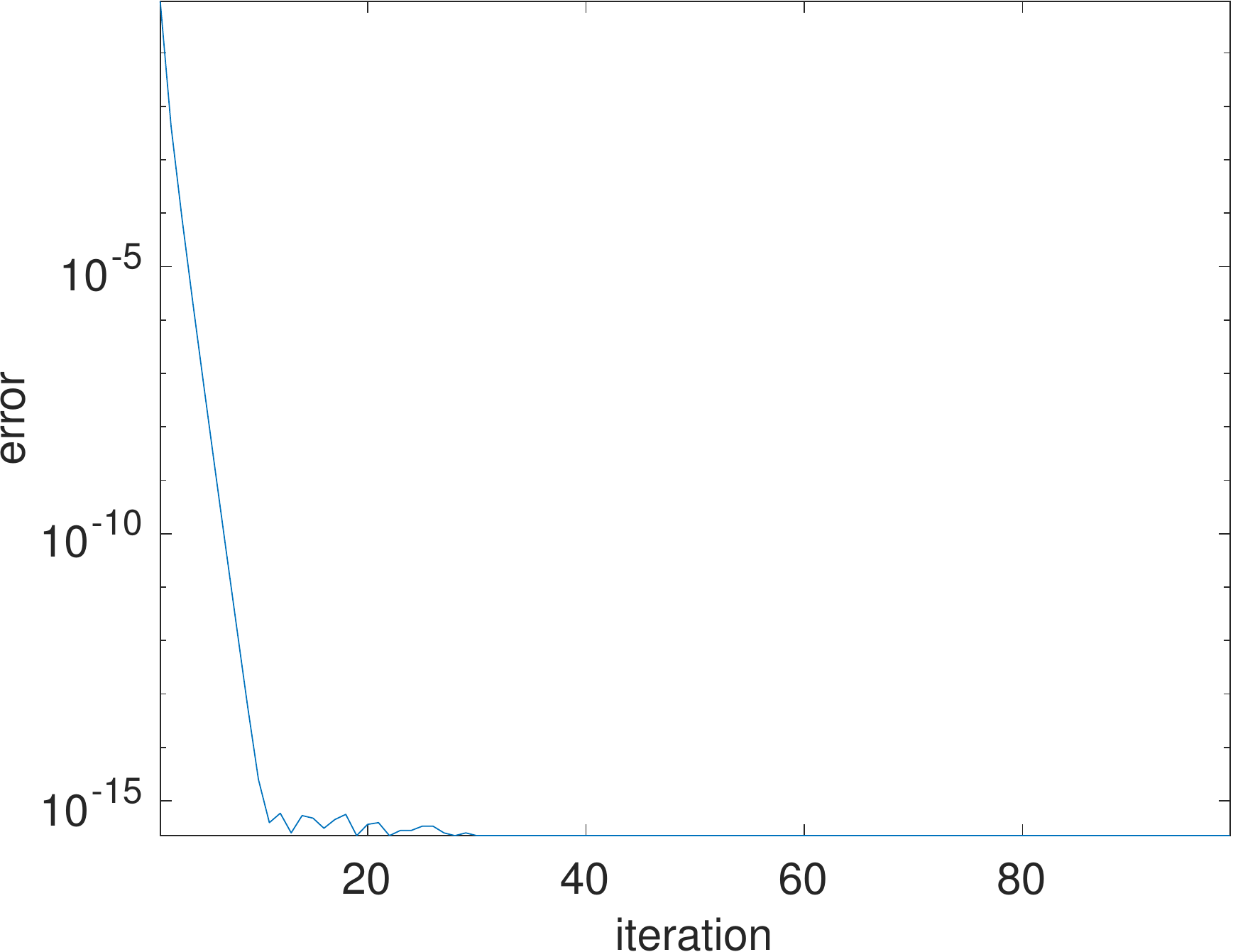}
  \includegraphics[width=0.32\textwidth]{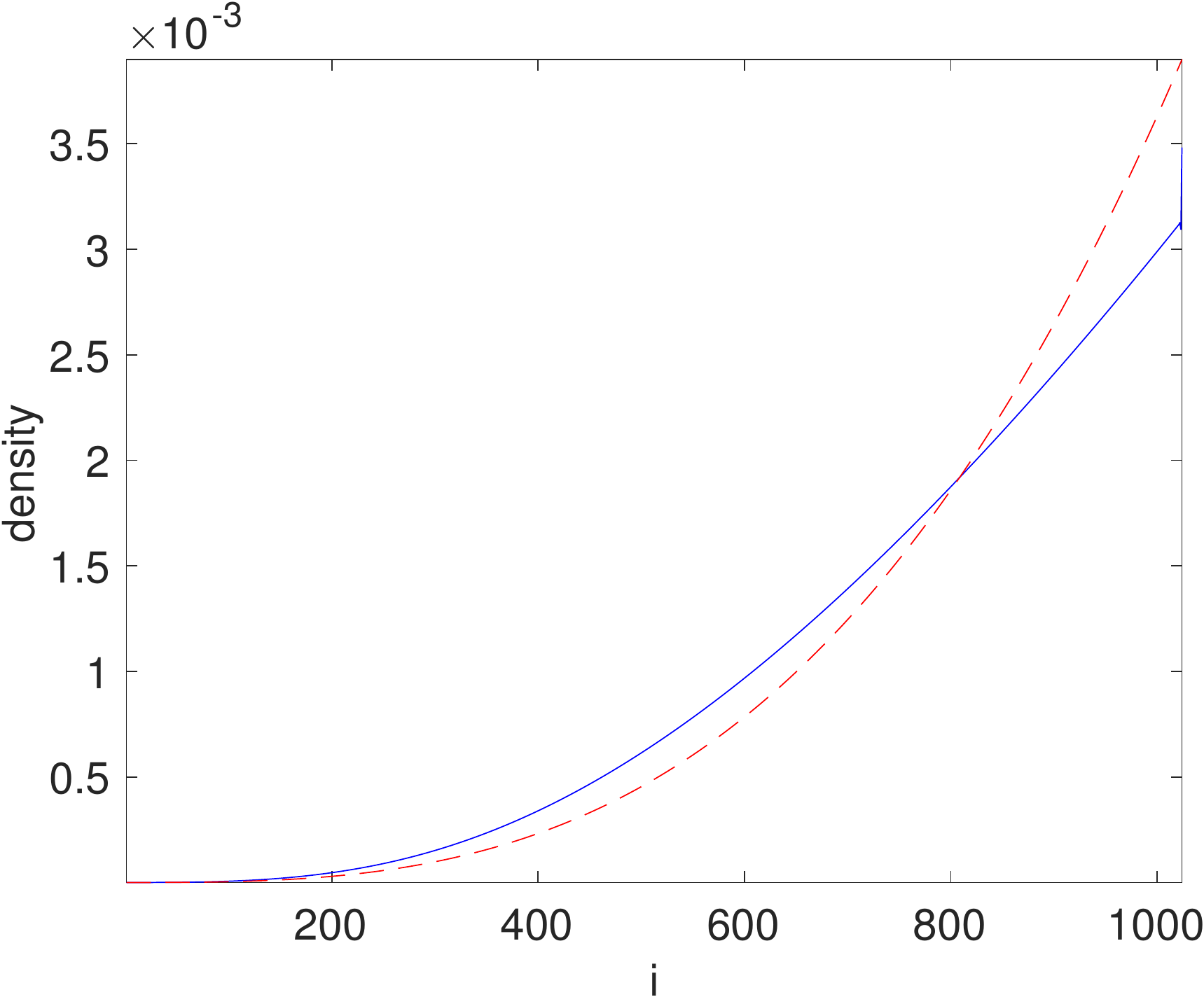}
  \caption{Reverse KL divergence, positive-definite case.  Left: free energy vs. iteration. Middle:
    free energy error vs. iteration. Right: density $p$ at the final iteration (solid line) and the
    minimizing density without $W$ term.}
  \label{fig:rkl_spd}
\end{figure}

As a numerical example, consider the periodic domain $[0,1]$ discretized with $n=1024$ points. The
potential $V$ is chosen to be zero and the interacting term is
\[
W_{ij} =
\begin{cases}
  \alpha, & i=j,\\
  \alpha/2, & i=j\pm 1,\\
  0, & \text{otherwise},
\end{cases}
\]
with $\alpha=10^2$. So $\alpha_i=10^2$ for each $i$. The step size $\Delta t$ is taken to be
$1$. Starting from a random initial condition, we run the descent algorithm for $100$ steps with the
results summarized in Figure \ref{fig:rkl_spd}. After about only $10$ iterations, the error is
reduced to about $10^{-15}$ .

%-----------------------------------
\section{Hellinger divergence} \label{sec:hlg}

For the Hellinger divergence
\[
D_{\Hlg}(p||\mu) = \sum_i (\sqrt{p_i}-\sqrt{\mu_i})^2 = -2 \sum_i \sqrt{\mu_i p_i} + \text{cst.}
\]
The free energy up to a constant is 
\[
F_{\Hlg}(p) = -2 \sum_i \sqrt{\mu_i p_i} + \sum_i V_i p_i + \frac{1}{2}\sum_{i,j} p_i W_{ij} p_j.
\]
The Hessian is given by
\[
\frac{\delta^2 F_{\Hlg}}{\delta p^2} = \diag\left( \frac{\mu^{1/2}}{2 p^{3/2}}\right) + W.
\]
Notice that the Hessian can be quite far from the mirror descent choice $\diag(1/(2p^{3/2}))$ even
when $W$ is zero, since $\mu$ can be drastically different for different $i$. When $W$ is
non-positive-definite, it is safe to continue using $\diag\left(\mu^{1/2}/(2p^{3/2})\right)$ as the
gradient metric. When $W$ is positive-definite, we extract the diagonal $\alpha = \diag(W)$ and use
$\diag\left( \mu^{1/2}/(2p^{3/2}) + \alpha\right)$ as the gradient metric.

%-------
\subsection{Non-positive-definite case}\label{sec:hlg_non}

Using $\diag\left(\mu^{1/2}/(2p^{3/2})\right)$ as the metric, the gradient flow is
\[
\dot{p} = -\frac{1}{\mu^{1/2}/(2p^{3/2})} \left(-\sqrt{\frac{\mu}{p}} + V + Wp +c \right).
\]
Moving the metric to the left hand side gives
\[
\dot{\left(-\sqrt{\mu/p}\right)} = -(-\sqrt{\mu/p} + V + Wp +c).
\]
If we introduce a reparameterization from $p\in\R^n$ to $g\in\R^n$ with
$g_i=\phi_i(p_i)\equiv-\sqrt{\mu_i/p_i}$ and $p_i=\phi_i^{-1}(g_i) = \mu_i/g_i^2$
\begin{align*}
  & \phi_i: p_i\rightarrow g_i,\quad (0,1)\rightarrow (-\infty, -\sqrt{\mu_i}),\\
  & \phi_i^{-1}: g_i\rightarrow p_i,\quad (-\infty,-\sqrt{\mu_i})\rightarrow (0,1),
\end{align*}
the gradient flow becomes
\[
\dot g = -(g + V + W p +c).
\]
The explicit Euler discretization gives
\begin{align*}
  &\tg = g^k -\Delta t (g^k + V + W p^k),\\
  & g^{k+1} = \tg+c.
\end{align*}
The constant $c$ is determined by the normalization condition
\[
\sum_i \phi_i^{-1}(\tg_i+c) = 1,
\]
which can be solved since it is monotone. The correct value $c$ can be shown to be in
\[
\left( \min\left(-\tg_i-\sqrt{n\mu_i}\right), \min(-\tg_i-\sqrt{\mu_i}) \right).
\]
Plugging the two endpoints of the interval shows that at the left endpoint
$\sum_i\phi_i^{-1}(\tg_i+c)<1$ and at the right endpoint $\sum_i\phi_i^{-1}(\tg_i+c)>1$. Therefore,
there is a unique $c$ value satisfies $\sum_i\phi_i^{-1}(\tg_i+c) = 1$ within this interval.

\begin{figure}[h!]
  \centering
  \includegraphics[width=0.32\textwidth]{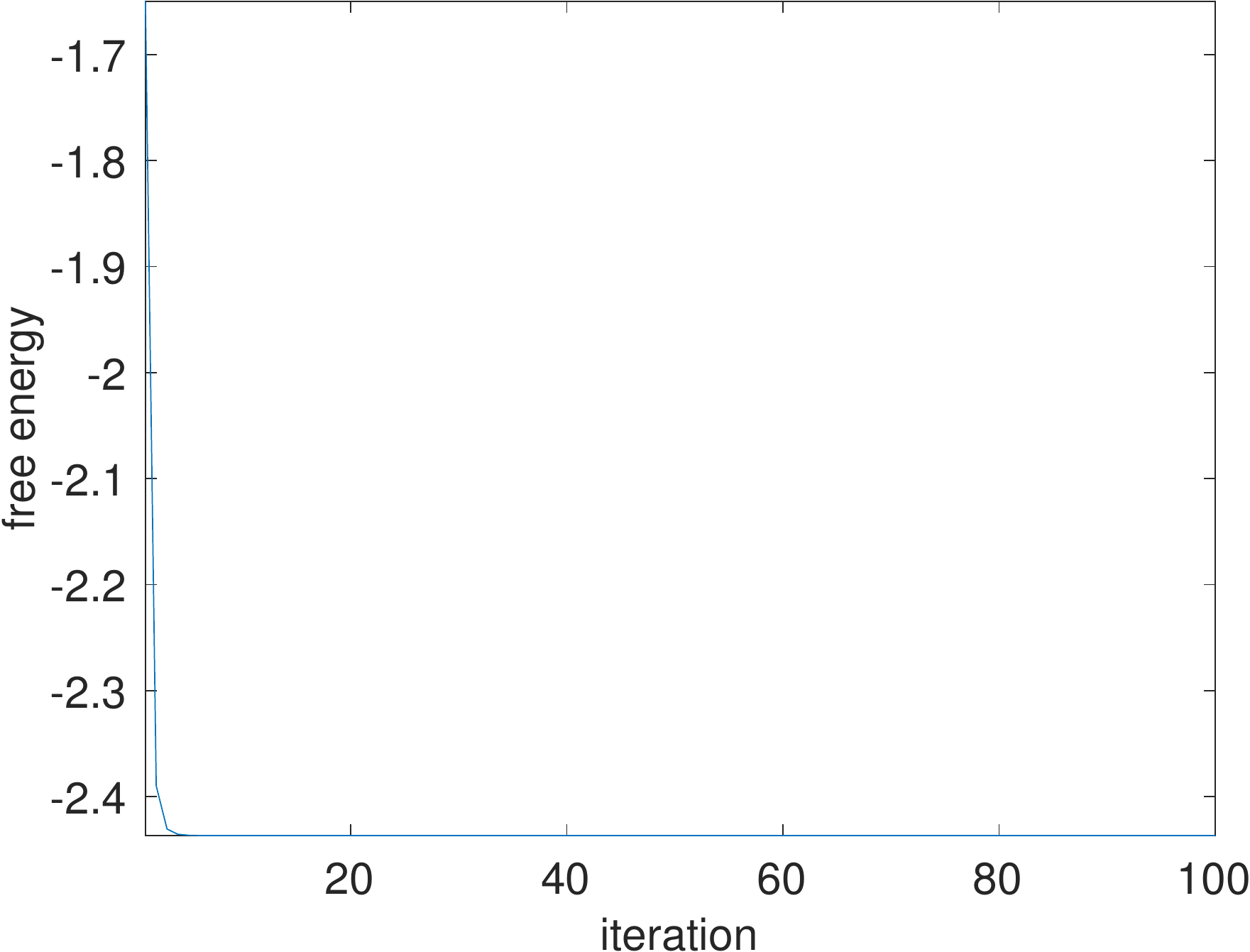}
  \includegraphics[width=0.32\textwidth]{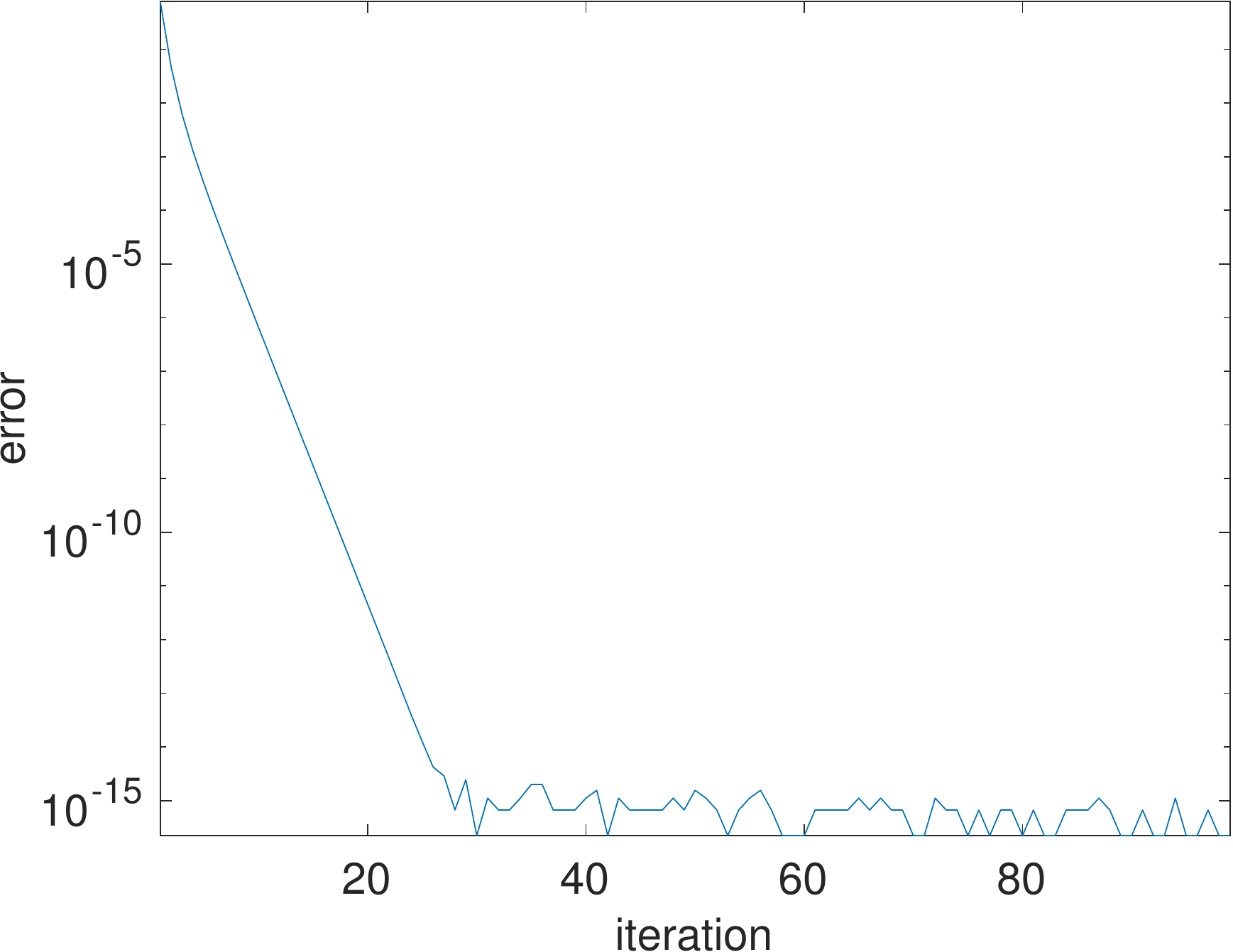}
  \includegraphics[width=0.32\textwidth]{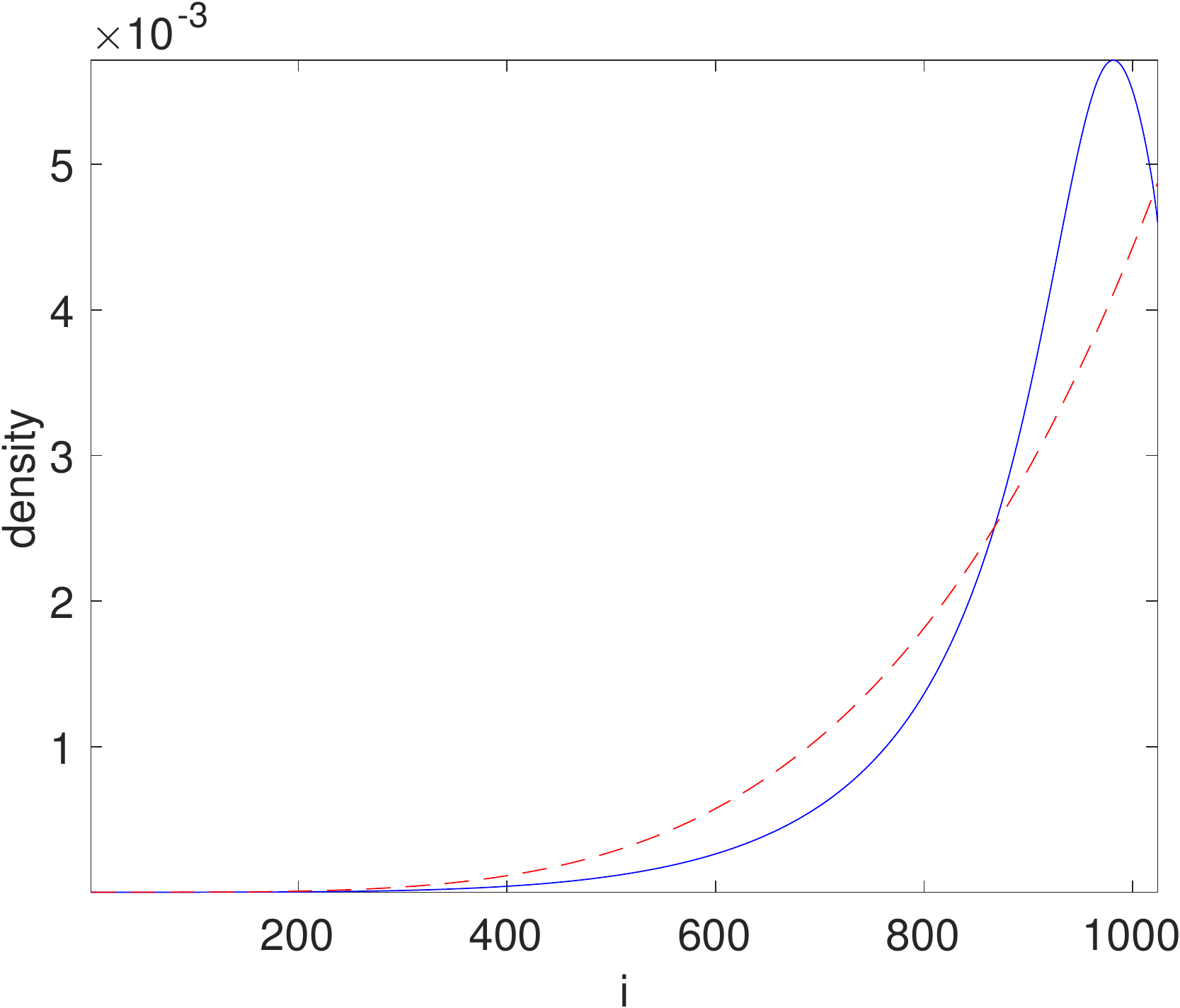}
  \caption{Hellinger divergence, non-positive-definite case with a Keller-Segel free energy.  Left:
    free energy vs. iteration. Middle: free energy error vs. iteration. Right: density $p$ at the
    final iteration (solid line) and the reference measure $\mu$ (dashed line). The reference
    density is the minimizer when $W=0$.}
  \label{fig:hlg_non}
\end{figure}

We illustrate the efficiency of the algorithm using a Keller-Segel model.  Consider the domain
$[0,1]$ discretized with $n=1024$ points $\{x_i=\frac{i}{n}\}$. The potential $V$ is zero and the
interacting term $W_{ij}$ is given by
\[
W_{ij} = \frac{1}{3} \ln(|x_i-x_j|+\varepsilon)
\]
with $\varepsilon=10^{-6}$. The reference measure $\mu_i$ is taken to be
\[
\mu_i \sim x_i^4.
\]
The step size $\Delta t$ is taken to be $1$. Starting from a random initial condition, we run the
descent algorithm for $100$ steps and the results are summarized in Figure \ref{fig:hlg_non}. Within
$30$ iterations, it reaches within $10^{-15}$ accuracy.

%-------
\subsection{Positive-definite case}\label{sec:hlg_spd}

Using $\diag\left( \mu^{1/2}/(2p^{3/2})+\alpha \right)$ as the metric, the gradient flow is
\[
\dot{p} = -\frac{1}{\mu^{1/2}/(2p^{3/2}) + \alpha} \left(-\sqrt{\frac{\mu}{p}} + V + Wp +c\right).
\]
Moving the metric to the left hand side gives
\[
\dot{\left(-\sqrt{\mu/p} +\alpha p\right)} = -(-\sqrt{\mu/p} + \alpha p + V + (W-\alpha)p +c).
\]
If we introduce a reparameterization from $p\in\R^n$ to $g\in\R^n$ with
$g_i=\phi_i(p_i)\equiv-\sqrt{\mu_i/p_i}+\alpha_i p_i$:
\begin{align*}
  & \phi_i: p_i\rightarrow g_i,\quad (0,1)\rightarrow (-\infty, -\sqrt{\mu_i}+\alpha_i),\\
  & \phi_i^{-1}: g_i\rightarrow p_i,\quad (-\infty,-\sqrt{\mu_i}+\alpha_i)\rightarrow (0,1),
\end{align*}
the gradient flow becomes
\[
\dot g = -(g + V + (W-\alpha) p +c).
\]
An explicit Euler discretization gives
\begin{align*}
  &\tg = g^k -\Delta t (g^k + V + (W-\alpha)p^k),\\
  & g^{k+1} = \tg+c.
\end{align*}
The constant $c$ is determined by the normalization condition
\[
\sum_i \phi_i^{-1}(\tg_i+c) = 1,
\]
which can be solved since it is monotone. The correct value $c$ can be shown to be in
\[
\left( \min\left(-\tg_i-\sqrt{n\mu_i}+\frac{\alpha_i}{n}\right), \min(-\tg_i-\sqrt{\mu_i}+\alpha_i) \right).
\]
Plugging the two endpoints of the interval shows that the left endpoint
$\sum_i\phi_i^{-1}(\tg_i+c)<1$ and at the right endpoint $\sum_i\phi_i^{-1}(\tg_i+c)>1$. Therefore,
there is a unique $c$ value satisfies $\sum_i \phi_i^{-1}(\tg_i+c) = 1$ within this interval. 

\begin{figure}[h!]
  \centering
  \includegraphics[width=0.32\textwidth]{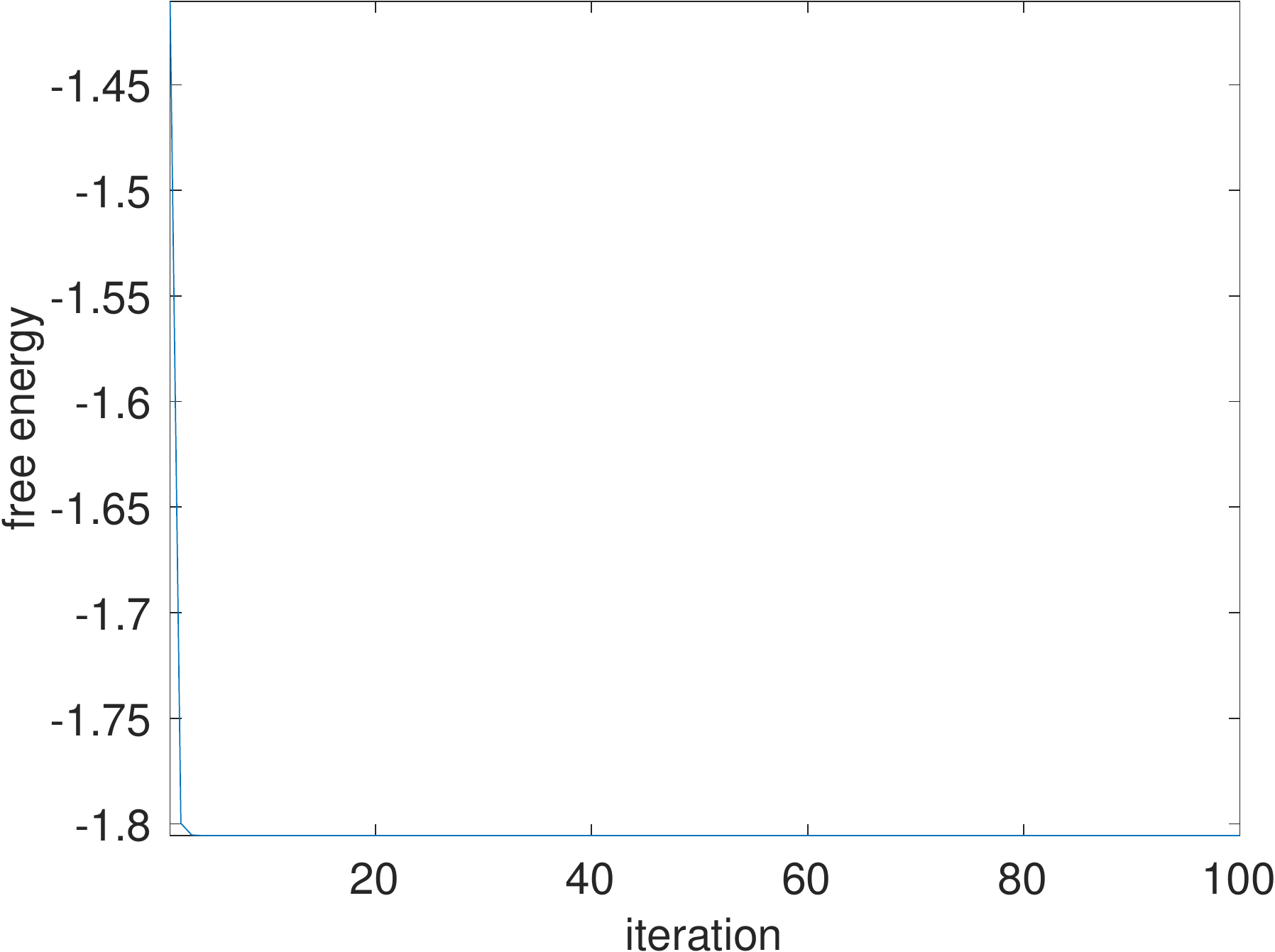}
  \includegraphics[width=0.32\textwidth]{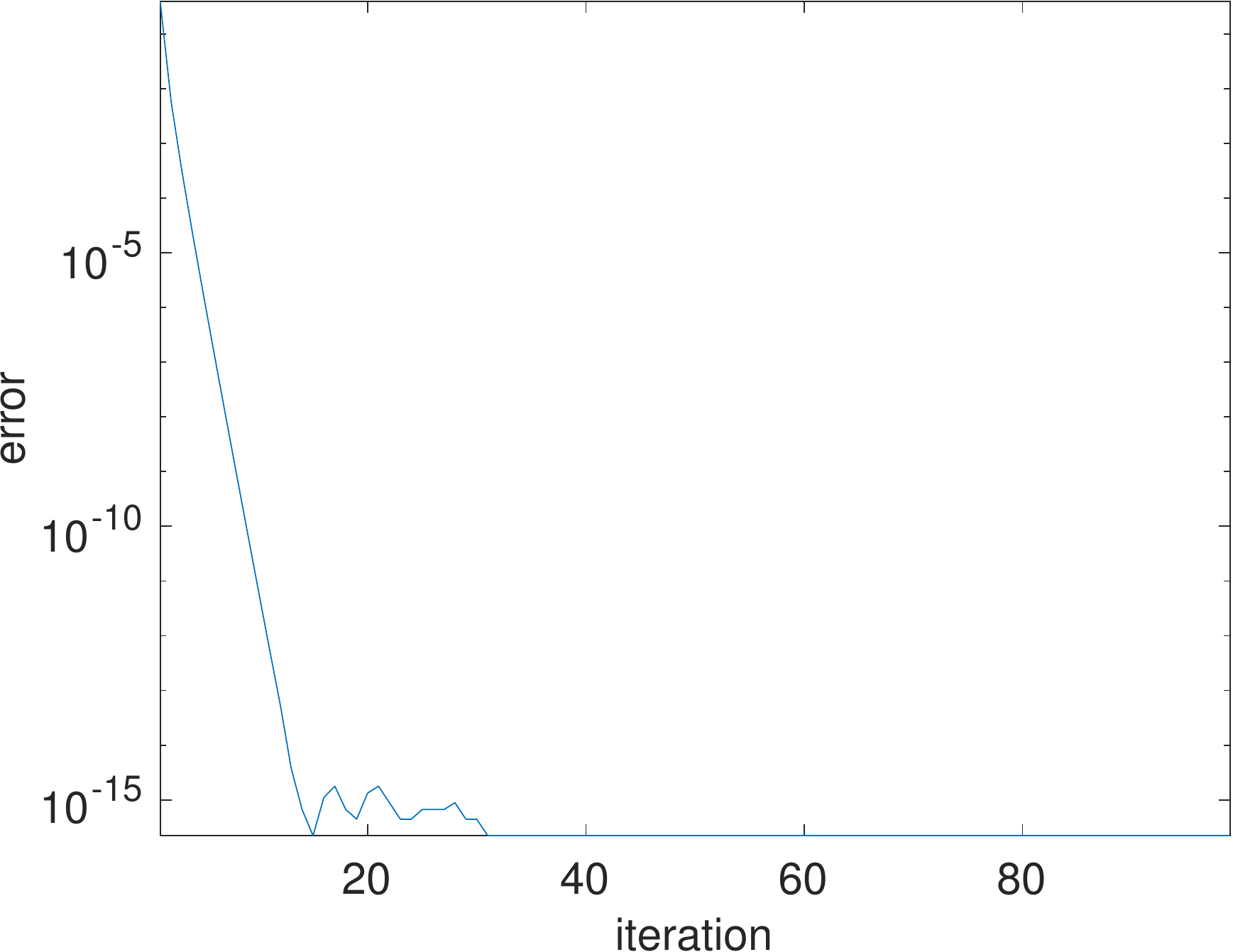}
  \includegraphics[width=0.32\textwidth]{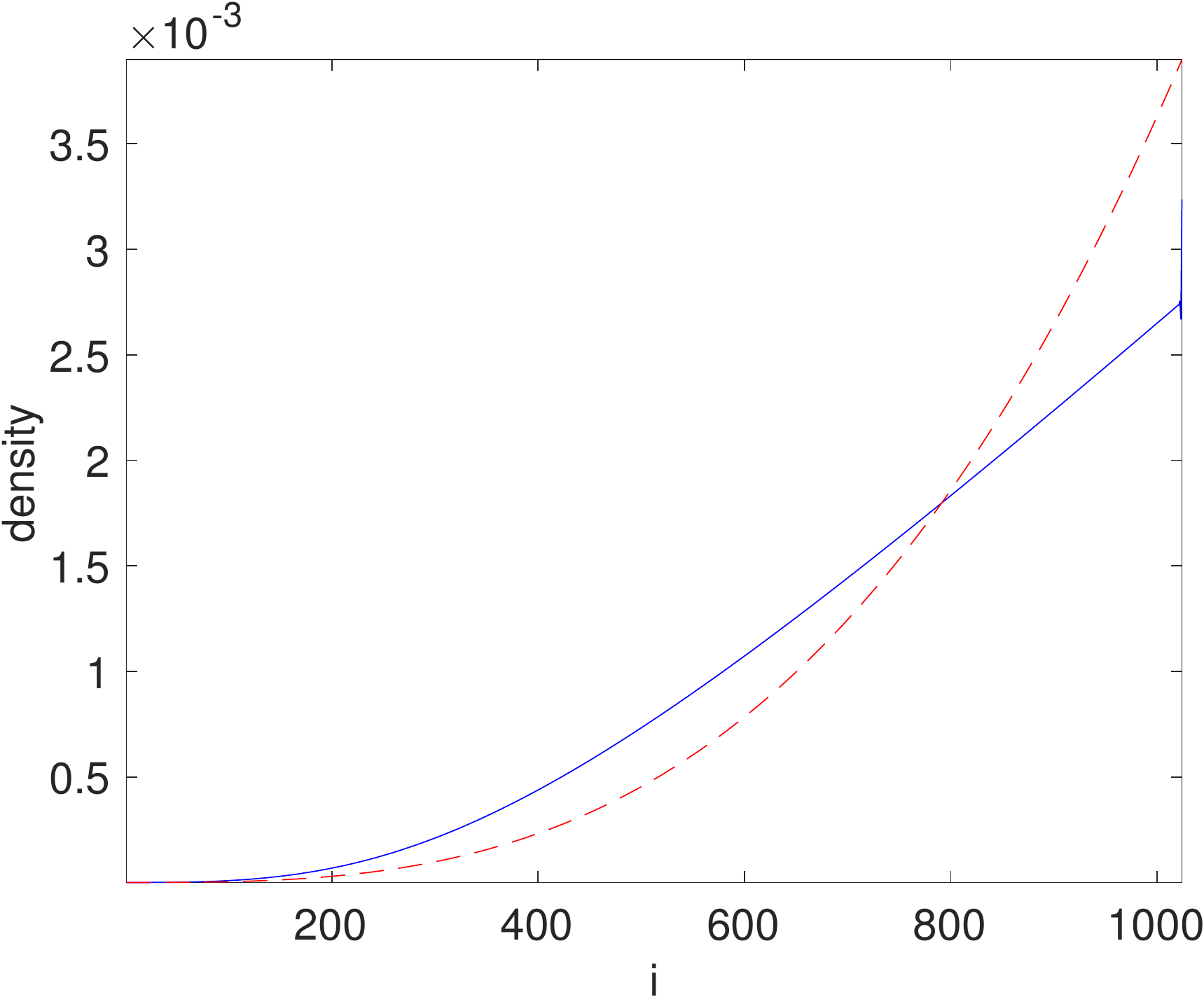}
  \caption{Hellinger divergence, positive-definite case.  Left: free energy vs. iteration. Middle:
    free energy error vs. iteration. Right: density $p$ at the final iteration (solid line) and the
    minimizing density without the $W$ term.}
  \label{fig:hlg_spd}
\end{figure}

In the numerical test, we consider the periodic domain $[0,1]$ discretized with $n=1024$ points. The
potential $V$ is chosen to be zero and the interacting term is
\[
  W_{ij} =
  \begin{cases}
    \alpha, & i=j,\\
    \alpha/2, & i=j\pm 1,\\
    0, & \text{otherwise},
  \end{cases}
\]
with $\alpha=10^2$. This leads to $\alpha_i=10^2$ for each $i=1,\ldots,n$. The step size $\Delta t$
is taken to be $1$. Starting from a random initial condition, we run the descent algorithm for $100$
steps. The results are summarized in Figure \ref{fig:hlg_spd}. Within about $15$ iterations, it
converges to an accuracy of order $10^{-15}$.

%-----------------------------------
\section{Discussions}

This paper proposes mirror-descent-type algorithms for minimizing interacting free energies. Below
we point out a few questions for future work. First, the proposed algorithms are obtained from
discretizing the continuous-time gradient flow with a new metric based on $\mu$ and $W$. One can
also derive the algorithm in a more traditional mirror descent form by starting from the
corresponding Bregman divergences.

Second, this paper considers three cases: KL divergence, reverse KL divergence, and Hellinger
divergence. In fact, the same procedure can be extended to most $\alpha$-divergences
\cite{amari2016information}.

When we treat the non-positive-definite case, $W$ is simply dropped in the design of the new
metric. A more accurate, but potentially more computational intensive, alternative is to find a
positive-definite approximation to $W$ and then combine it with the Hessian from the divergence
term.

This interacting term of the free energy considered in this paper is only of quadratic form. It is
plausible that a similar procedure can be developed for non-quadratic interacting terms, as long as
there is an efficient way to approximate the diagonal of the Hessian.

\bibliographystyle{abbrv}

\bibliography{ref}

\end{document}